\documentclass[psamsfont,a4paper,11pt]{amsart}

\usepackage[english]{babel}   %loads english

\usepackage{setspace}     %controls line-spacing
\usepackage{enumerate} 	%allows lists
\usepackage{afterpage} 	%keeps tables and figures where they belong
\usepackage[usenames,dvipsnames]{color}
\usepackage{graphicx}	%allows inclusion of pictures

\usepackage{amsmath} 	%just apparently the most important thing ever
\usepackage{amssymb}	%symbols in math
\usepackage{mathrsfs} 	%really fancy script font
\usepackage{amsfonts} 	%standard fonts
\usepackage{latexsym} 	%more symbols
\usepackage[latin1]{inputenc} %sets up proper special characters

\usepackage{amsthm} 	%sets up theorem styling
\usepackage{stackrel}       %allows stacking in mathmode  
\usepackage{amscd} 	%commutative diagrams                                 

\usepackage{tabularx}	%automatically widening tables
\usepackage{longtable}	%for tables longer than a page

\usepackage{verbatim} %comment environment?
\usepackage{blindtext} %something about dummy text

\usepackage{lipsum}

\newcommand\blfootnote[1]{%
  \begingroup
  \renewcommand\thefootnote{}\footnote{#1}%
  \addtocounter{footnote}{-1}%
  \endgroup
}

\allowdisplaybreaks

\def\cal{\mathcal}

\def\O{{\cal O}}

\newcommand{\R}{\mathbb{R}} 
 
\newcommand{\NN}{\mathbb{N}} 
\newcommand{\Z}{\mathbb{Z}} 
\newcommand{\Q}{\mathbb{Q}}

\newtheorem{Thm}{Theorem}[section]		%%with numbering
\newtheorem{Lemma}[Thm]{Lemma} 
 
\newtheorem{Prop}[Thm]{Proposition}

\newtheorem*{cor}{Corollary}

\theoremstyle{definition}

\theoremstyle{remark}
 
\newtheorem*{rmk}{Remark}

\newtheorem*{ex}{Example}

\newtheorem{ind}[]{{\rm\it Indice}}

\newcommand{\overbar}[1]{\mkern 1.5mu\overline{\mkern-1.5mu#1\mkern-1.5mu}\mkern 1.5mu}

\usepackage{multicol}

\usepackage{tikz}
\usetikzlibrary{arrows}

\usepackage{paralist}
\usepackage{cite}

\title{Conjugacy growth series for wreath product finitary symmetric groups}

\author[Wagner]{Ian Wagner}

\begin{document}

\maketitle

\blfootnote{$2010$ \textit{Mathematics Subject Classification}.  $11$F$33$, $11$F$37$, $20$B$30$, $20$C$32$. }

\begin{abstract}
In recent work, Bacher and de la Harpe define and study conjugacy growth series for finitary permutation groups.  In two subsequent papers, Cotron, Dicks, and Fleming study the congruence properties of some of these series.  We define a new family of conjugacy growth series for the finitary alternating wreath product that are related to sums of modular forms of integer and half-integral weights, the so-called \textit{mixed weight modular forms}.  The previous works motivate the study of congruences for these series.  We prove that congruences exist modulo powers of all primes $p \geq 5$.  Furthermore, we lay out a method for studying congruence properties for sums of mixed weight modular forms in general.  

\end{abstract}

\section{Introduction and Statement of Results}

In the recent paper \cite{BD}, Roland Bacher and Pierre de la Harpe develop the theory of conjugacy growth series.  This theory uses the minimum word length statistics and the conjugacy classes of a group to produce the conjugacy growth series.  In particular, this theory ties together infinite permutation groups with finite support and the usual number theoretic partition function.

Let $G$ be a group and $S$ a set that generates $G$, then for each $g \in G$ define the \textit{word length}, $\ell_{G,S}(g)$, to be the smallest nonnegative integer $n$ for which there are $s_{1}, s_{2}, ..., s_{n} \in S \cup S^{-1}$ such that $g = s_{1}s_{2}\cdot \cdot \cdot s_{n}$.  Define the \textit{conjugacy length}, $\kappa_{G,S}(g)$, as the smallest integer $n$ for which there exists $h$ in the conjugacy class of $g$ such that $\ell_{G,S}(h) = n$.  For $n \in \NN$ define $\gamma_{G,S}(n) \in \NN \cup \{0\} \cup \{\infty\}$ as the number of conjugacy classes of $G$ which contain elements $g$ with $\kappa_{G,S}(g) = n$.  Whenever $\gamma_{G,S}(n)$ is finite for all $n$, we can define the \textit{conjugacy growth series}:
\[C_{G,S}(q) = \sum_{n=0}^\infty \gamma_{G,S}(n)q^{n} = \sum_{g \in Conj(G)} q^{\kappa_{G,S}(g)} \in \NN[[q]], \tag{1.1} \]
where the second sum is taken over representatives of conjugacy classes of $G$. 

We call Sym($\NN$) the \textit{finitary symmetric group} of $\NN$.  It is the group of permutations of $\NN$ with finite support.  Let the \textit{finitary alternating group} of $\NN$, Alt($\NN$), be the subgroup of Sym($\NN$) of permutations with even signature.  Define the two generating sets of Sym($\NN$), $S_{\NN}^{Cox} = \{(i, i+1) : i \in \NN \}$ and $T_{\NN} = \{(x,y) : x, y \in {\NN}$ are distinct$\}$.  Let $S \subset$ Sym($\NN$) be a generating set such that $S_{\NN}^{Cox} \subset S \subset T_{\NN}$.  Then Bacher and de la Harpe prove that (see Proposition $1$ of \cite{BD})
\[ C_{\rm{Sym}(\NN),S}(q) = \sum_{n=0}^\infty p(n)q^n = \prod_{n=1}^\infty \frac{1}{1-q^{n}}. \tag{1.3} \]

Similarly we can define two generating sets of Alt($\NN$), $S_{\NN}^{A} = \{(i, i+1, i+2) : i \in \NN \}$ and $T_{\NN}^{A} = \bigcup_{g \in \rm{Alt}(\NN)} gS_{\NN}^{A}g^{-1}$, where $T_{\NN}^{A}$ is the subset of all 3-cycles.

Let $S^{'} \subset$ Alt($\NN$) be a generating set such that $S_{\NN}^{A} \subset S^{'} \subset T_{\NN}^{A}$.  Then  Bacher and de la Harpe also prove that (see Proposition $11$ of \cite{BD})
\[C_{\rm{Alt}(\NN),S^{'}}(q) = {\sum_{n=0}^\infty p(n)q^{n}}{\sum_{m=0}^\infty p_{e}(m)q^{m}} = \frac{1}{2} \prod_{n=1}^\infty \frac{1}{(1-q^{n})^{2}} + \frac{1}{2} \prod_{n=1}^\infty \frac{1}{1-q^{2n}},  \tag{1.4}\]
where $p_{e}(m)$ denotes the number of partitions of $m$ into an even number of parts. 

\begin{comment}
This last line can also be viewed in terms of eta-products:
\[\frac{1}{2} \prod_{n=1}^\infty \frac{1}{(1-q^{n})^{2}} + \frac{1}{2} \prod_{n=1}^\infty \frac{1}{1-q^{2n}}\]
\[=\frac{1}{2} \frac{q^{1/12}}{\eta(z)^{2}} + \frac{1}{2} \frac{q^{1/12}}{\eta(2z)} \]
where $\eta(z) = q^{1/24} \prod_{n=1}^\infty (1-q^{n})$.  This is of particular interest because, up to $q^{1/12}$, it is a sum of a modular form of weight $-1$ and a modular form of weight $-\frac{1}{2}$.  In \cite{CDF2} Cotron, Dicks, and Fleming prove asymptotics and congruences for this sum.
\end{comment}

In Section $2$ we will generate conjugacy growth series, $C_{W_{M}^{'},S_{*}^{'}}(q)$, that are powers of equation $(1.4)$.  If $M$ is a positive integer, let \[ \sum_{n=0}^\infty \gamma_{W_{M}^{'}, S_{*}^{'}}(n) q^{n} = \left( \frac{1}{2} \prod_{n=1}^\infty \frac{1}{(1-q^{n})^{2}} + \frac{1}{2} \prod_{n=1}^\infty \frac{1}{1-q^{2n}} \right)^{M}.  \tag{1.5}\]  
These generating functions turn out to be the conjugacy growth series for wreath products of Alt($\NN$).

Let $H$ be a group.  Then $W := H \wr_{\NN}$Sym($\NN$) $=H^{(\NN)} \rtimes$ Sym($\NN$) is called a \textit{permutation wreath product}.  Denote by $H^{(\NN)}$ the group of functions from $\NN$ to $H$.   Sym($\NN$) has a natural action on $H^{(\\N)}$; $\sigma \in$ Sym($\NN$) acts on $\phi \in H^{(\NN)}$ by $\sigma(\phi) = \phi \circ \sigma^{-1}$.  One can also think of $H^{(\NN)}$ as $|\NN|$ copies of $H$, and so an element of $H^{(\NN)}$ can be thought of as $|\NN|$ elements of $H$ indexed by $\NN$.  In particular, Sym($\NN$) acts naturally on these indices.  For $\sigma, \tau \in$ Sym($\NN$) and $\phi, \psi \in H^{(\NN)}$, the multiplication in the wreath product is given by $(\phi, \sigma)(\psi, \tau) = (\phi \sigma(\psi), \sigma \tau)$.  The alternating wreath product, $W^{'} := H \wr_{\NN}$Alt($\NN$), can be defined analogously.

\begin{comment}
The following theorems provide asymptotics and congruences for these series.   

\begin{Thm} \label{five}
If $M$ is a positive integer, let \[ \sum_{n=0}^\infty c_{M}(n) q^{n} = \Big( \frac{1}{2} \prod_{n=1}^\infty \frac{1}{(1-q^{n})^{2}} + \frac{1}{2} \prod_{n=1}^\infty \frac{1}{1-q^{2n}} \Big)^{M}.  \tag{1.6}\]
Then, as $n \rightarrow \infty$, we have that \[c_{M}(n) \sim \frac{M^{\frac{2M+1}{4}}}{2^{2M+1} 3^{\frac{2M+1}{4}} n^{\frac{2M+3}{4}}} e^{2 \pi \sqrt{\frac{Mn}{3}}}. \]

\end{Thm}
\end{comment}

In view of $(1.5)$ and its interpretation in terms of conjugacy growth series, it is natural to study the congruence properties of the coefficients of these functions in the spirit of the earlier work of Bacher and de la Harpe and Cotron, Dicks, and Fleming in \cite{BD} and \cite{CDF2}. In \cite{CDF2} Cotron, Dicks, and Fleming only discuss congruences for $C_{W_{M}^{'},S_{*}^{'}}(q)$ for $M=1$ and for powers of the primes $5$ and $7$.  For example, they proved that 
\[\gamma_{W_{1}^{'},S_{*}^{'}}(2 \cdot 5^4 n + 1198) \equiv 0 \pmod{5}\] and 
\[\gamma_{W_{1}^{'},S_{*}^{'}}(2 \cdot 7^6 n + 225494) \equiv 0 \pmod{49}.\]  It is natural to ask for a more complete description of congruences for all of the general wreath products.  This also motivates the study of sums of mixed weight modular forms in general.  

Cotron, Dicks, and Fleming use the theory of modular forms to obtain their results.  Therefore, one expects to use modular forms to study the conjugacy growth series $(1.5)$.  However, a difficulty arises; these functions are \textit{mixed weight modular forms}, finite sums of modular forms with different weights.  Therefore, we must first obtain a general theorem about congruences for coefficients of mixed weight modular forms.  Let $q := e^{2 \pi i z}$ with $z \in \mathbb{H}$, and let $\mathcal{M}_{k}(\Gamma_{0}(N), \chi)$ denote the space of weakly holomorphic modular forms of weight $k$, level $N$, and with character $\chi$.  \textit{Weakly holomorphic modular forms} are those meromorphic modular forms whose poles (if any) are supported at cusps.  Then the following offers such a theorem.

\begin{Thm} \label{main}

Let $K$ be an algebraic number field with ring of integers $\mathcal{O}_{K}$.  Suppose $f_{i}(z) = \sum_{n=0}^\infty a_{i}(n) q^{n} \in \mathcal{M}_{k_{i}}(\Gamma_{0}(N_{i}), \chi_{i}) \cap \mathcal{O}_{K}((q))$, $g_{j}(z) = \sum_{m=0}^\infty b_{j}(m) q^{m} \in \mathcal{M}_{\lambda_{j} + \frac{1}{2}}(\widetilde{\Gamma_{0}(M_{j})}, \chi_{j}) \cap \mathcal{O}_{K}((q))$ where $4 \mid M_{j}$ for every $j$, and let
\[F(z) = \sum A(n)q^n = \sum_{i=1}^{u} f_{i}(z) + \sum_{j=1}^{v} g_{j}(z). \]
Let $N$ be minimal such that $N_{i} \mid N$ and $M_{j} \mid N$ for every $i$ and $j$, and let $p$ be prime such that $(N, p) = 1$.  If $r$ is a sufficiently large integer, then for each positive integer $j$, a positive proportion of primes $Q \equiv -1 \pmod{N p^j}$ have the property that
\[A(Q^{4t+3} p^r n) \equiv 0 \pmod{p^j},\]
where $(Qp, n) = 1$ and $t$ is a nonnegative integer. 

\end{Thm}

\begin{rmk}
If there are any half integral weight forms in the sum of forms above, then we will have $4 \vert N$.  In this case, it is clear that $p$ must be an odd prime.
\end{rmk}

Applying Theroem \ref{main} to the conjugacy growth series leads to the following theorem.

\begin{Thm} \label{six}
Suppose $p \geq 5$ is prime and let $j$ be a positive integer.  If $r$ is a sufficiently large integer, then for a positive proportion of primes $Q \equiv -1 \pmod{576p^{j}}$, we have that 
\[\gamma_{W_{M}^{'},S_{*}^{'}} \left(\frac{Q^{4t+3}p^{r}n + M}{12} \right) \equiv 0 \pmod{p^{j}}\]
for all $n$ coprime to $Qp$, and for all nonnegative integers $t$.

\end{Thm} 

A technique from \cite{Ono2} can be adapted in order to do a computer search for examples of congruences.  Using the work in \cite{Ono2} we find that  
\[\gamma_{W_{2}^{'}, S_{*}^{'}}\left(\frac{7 n + 2}{12}\right) \equiv \eta^{10}(24z) \pmod{7},\]
where $\eta(z)$ is \textit{Dedekind's eta-function}.  From this we find the following congruence example.

\begin{ex}
We have that
\[\gamma_{W_{2}^{'}, S_{*}^{'}}\left(\frac{7 n + 2}{12}\right) \equiv 0 \pmod{7}\]
whenever $n \not\equiv 10 \pmod{24}$.  Moreover, the above congruence is true when $ n = 24t + 10$ and $t \equiv 2, 4, 5,$ or $6 \pmod{7}$.

\end{ex}

To prove Theorems \ref{main} and \ref{six} we will make use of the theory of modular forms.  The relevant generating functions  for Theorem \ref{six} turn out to be mixed weight modular forms.  We will use the work of Treneer in \cite{Tre} on weakly holomorphic modular forms and a famous theorem of Serre (see \cite{ONO1}).  We will show that Theorem \ref{main} follows from a proposition of Ono and Skinner in \cite{OS} which allows us to use the theory of Galois representations attached to modular forms for a finite set of modular forms simultaneously.  Section $3.1$ will cover basic facts about modular forms, operators for modular forms, and cusps.  In Section $3.2$ we will state important propositions of Treneer, Serre, Ono and Ahlgren, and others which are vital to our proofs.  We will also discuss the interplay between Galois representations and Hecke operators in Section $3.2$.  Theorem \ref{main} will be proved in Section $4$ and Theorem \ref{six} in will be proved Section $5$.  The paper will conclude with a short explanation of the example at the end of Section $1$ in Section $6$.

\section*{Acknowledgements}
The author would like to thank Ken Ono for his guidance on this project and Maddie Locus for her helpful commments.

\section{Conjugacy growth series for the finitary alternating wreath product}

It is natural to ask if there are other subgroups of the finitary symmetric group that produce interesting conjugacy growth series.  Recall the wreath product $W = H \wr_{\NN}$Sym($\NN$). \begin{comment} $W$ acts faithfully on $H \times \NN$ by $(\phi, \sigma)$ acting by $(h, n) \mapsto (\phi(\sigma(n))h, \sigma(n))$. \end{comment}

For $a \in H\setminus\{1\}$ and $m \in \NN$, let $\phi_{m}^{a} \in W$ be the permutation that maps $(h, m) \in H \times \NN$ to $(ah, m)$ and fixes $(h, n)$ if $n \neq m$. Note that $(\phi_{m}^{a})_{a \in H\setminus \{1\}, m \in \NN}$ generates the subgroup $H^{(\NN)}$ and that $\phi_{m}^{a}$ and $\phi_{k}^{b}$ are conjugate in $W$ if and only if $a$ and $b$ are conjugate in $H$.  For $m \in \NN$, let $H_{m} = \{\phi_{m}^{a} : a \in H\setminus\{1\}\}$ and let $T_{H} = \bigcup_{m \in \NN} H_{m}$ be a subset of $H^{(\NN)}$.  Recall that $T_{\NN}$ is the subset of all transpositions in Sym($\NN$).  We consider subsets $S_{H} \subset T_{H}$ and $S_{\NN} \subset T_{\NN}$ and define $S_{*}$ to be the disjoint union $S_{H} \sqcup S_{\NN}$.  If $S_{H} = \{\phi_{m_{1}}^{a_{1}}, ..., \phi_{m_{r}}^{a_{r}}\}$ where $\{a_{1}, ..., a_{r}\}$ generate $H$, then $S_{*}$ generates $W$.  Define $S_{*}^{'}$ analogously using subsets of $T_{H}^{A}$ and $T_{\NN}^{A}$; then $S_{*}^{'}$ generates $W^{'} = H\wr_{\NN}$\rm{Alt}($\NN$).  This leads to the following proposition.

\begin{comment}

\begin{Prop}[Bacher-de la Harpe] \label{three}
Let $H$ be a finite group; denote by $M$ the number of conjugacy classes of $H$.  If $W = H\wr_{\NN}$\emph{Sym}($\NN$) and $S_{*}$ is a generating set satisfying (PCwr) in \cite{BD}, then
\[C_{W,S_{*}}(q) = \prod_{k=1}^\infty \frac{1}{(1-q^{k})^{M}}. \tag{1.5} \]

\end{Prop}

\end{comment}

\begin{Prop} \label{four}
Let $H$ be a finite group; denote by $M$ the number of conjugacy classes of $H$.  If $W_{M}^{'} = H\wr_{\NN}$\emph{Alt}($\NN$) and $S_{*}^{'}$ is a generating set satisfying (PCwr) in \cite{BD}, then
\[C_{W_{M}^{'},S_{*}^{'}}(q) = \left( \frac{1}{2} \prod_{n=1}^\infty \frac{1}{(1-q^{n})^{2}} + \frac{1}{2} \prod_{n=1}^\infty \frac{1}{1-q^{2n}} \right)^{M}. \]

\end{Prop}

\emph{Proof of Proposition \ref{four}:} For each $w' = (\phi, \sigma) \in W_{M}^{'} = H \wr_{\NN} \rm{Alt}(\NN)$ we can split $\sigma$ into the product of an even number of cycles of even length, $\sigma_{e}$, and the product of cycles of odd length, $\sigma_{o}$, so that $w' = (\phi, \sigma_{e} \sigma_{o})$.  We can associate each conjugacy class in $W_{M}^{'}$ to an $H_{*}$-indexed family of partitions.  Using the same notation as in \cite{BD} we associate the conjugacy classes in $H$ to the family of partitions 
\[\Big(\lambda^{(1)}, \nu^{(1)}; (\mu^{(\eta)}, \gamma^{(\eta)})_{\eta \in H_{*} \setminus 1} \Big),\]
where $\nu^{(1)}$ and $\gamma^{(\eta)}$ each have an even number of positive parts, in the following way.

Let $\NN^{(w)}$ be the finite subset of $\NN$ that is the union of the supports of $\phi$ and $\sigma$ and let $\sigma$ be the product of the disjoint cycles $c_{1}, ..., c_{k}$  where $c_{i} = (x_{1}^{(i)}, x_{2}^{(i)}, ..., x_{v_{i}}^{(i)})$ with $x_{j}^{(i)} \in \NN^{(w)}$ and $v_{i} = \text{ length}(c_{i})$.  We include cycles of length $1$ for each $n \in \NN$ such that $n \in $ sup$(\phi)$ and $n \notin $ sup$(\sigma)$ so that $\NN^{(w)} = \sqcup_{1\leq i \leq k}$ sup$(c_{i})$.  Define $\eta_{*}^{w}(c_{i}) \in H_{*}$ to be the conjugacy class of $\phi(x_{v_{i}}^{(i)}) \phi(x_{v{i}-1}^{(i)}) \cdot \cdot \cdot \phi(x_{1}^{(i)}) \in H$.  For $\eta \in H_{*}$ and $\ell \geq 1$, let $m_{\ell}^{w,\eta}$ denote the number of cycles $c$ in $\{c_{1}, ..., c_{k} \}$ that are of length $\ell$ and such that $\eta_{*}^{w}(c) = \eta$.  Let $\mu^{w, \eta} \vdash n^{w, \eta}$ be the partition with $m_{\ell}^{w, \eta}$ parts equal to $\ell$, for all $\ell \geq 1$.  Note that $\sum_{\eta \in H_{*}, \ell \geq 1} n^{w, \eta} = \sum_{\eta \in H_{*}, \ell \geq 1} \ell m_{\ell}^{w, \eta} = |\NN^{(w)}|$.  Also observe that the partition $\mu^{w, 1}$ does not have parts of size $1$ because if $v_{i} = 1$ then $\eta_{*}^{w}(c_{i}) \neq 1$.  Using the same notation as above, let $\lambda^{w, 1}$ be the partition with $m_{\ell}^{w, 1}$ parts equal to $\ell - 1$.  Because we are working in Alt($\NN$), we can write $\sigma = \sigma_{e} \sigma_{o}$; and so this method actually splits to map to two partitions, one of which has an even number of parts.  Define the \textit{type} of $w$ as the family $\Big(\lambda^{(1)}, \nu^{(1)}; (\mu^{(\eta)}, \gamma^{(\eta)})_{\eta \in H_{*} \setminus 1} \Big)$.  Then two elements in $W_{M}^{'}$ are conjugate if and only if they have the same type.  Thus, each $H_{*}$-indexed family of partitions, $\Big(\lambda^{(1)}, \nu^{(1)}; (\mu^{(\eta)}, \gamma^{(\eta)})_{\eta \in H_{*} \setminus 1} \Big)$,  is the type of  one conjugacy class in $W_{M}^{'}$.

Consider an $H_{*}$-indexed family of partitions $\Big(\lambda^{(1)}, \nu^{(1)}; (\mu^{(\eta)}, \gamma^{(\eta)})_{\eta \in H_{*} \setminus 1} \Big)$ and the corresponding conjugacy class in $W_{M}^{'}$.  Let $u^{(1)}, v^{(1)}, u^{(\eta)}, v^{(\eta)}$ be the sums of the parts of $\lambda^{(1)}, \nu^{(1)}, \mu^{(\eta)}, \gamma^{(\eta)}$ and let $k^{(1)}, t^{(1)}, k^{(\eta)}, t^{(\eta)}$ be the be the number of parts of  $\lambda^{(1)}, \nu^{(1)}, \mu^{(\eta)}, \gamma^{(\eta)}$ respectively.

Choose a representative $w' = (\phi, \sigma)$ of this conjugacy class such that \[\sigma = \prod_{i=1}^{k} c_{i} = \prod_{i=1}^{k} (x_{1}^{(i)}, x_{2}^{(i)}, ..., x_{\mu_{i}}^{(i)})\] and 
\[ \phi(x_{j}^{(i)}) = 1 \in H \text{ for all } j \in \{1, ..., \mu_{i} \}  \quad \text{ when }  \eta_{*}^{w}(c_{i}) = 1 \]
\[ \phi(x_{j}^{(i)}) = \begin{cases}
1 \text{ for all } j \in \{1, ..., \mu_{i} - 1\} \\
h \neq 1 \text{ for } j = \mu_{i} 
\end{cases}  \text{ when } \eta_{*}^{w}(c_{i}) \neq 1. \]
Observe that 
\[k = k^{(1)} + t^{(1)} + \sum_{\eta \in H_{*} \setminus 1} (k^{(\eta)} + t^{(\eta)} )\]
\[|\NN^{(w')}| = u^{(1)} + k^{(1)} + v^{(1)} + t^{(1)} + \sum_{\eta \in H_{*} \setminus 1} (u^{(\eta)} + v^{(\eta)} ).\]

\begin{comment} NEED MORE DETAILS
\end{comment}

Hence, the contribution  to $C_{W_{M}^{'},S_{*}}(q)$ from $\Big(\lambda^{(1)}, \nu^{(1)}; (\mu^{(\eta)}, \gamma^{(\eta)})_{\eta \in H_{*} \setminus 1} \Big)$ is \[\Big(q^{u^{(1)}}q^{v^{(1)}} \prod_{\eta \in H_{*} \setminus 1} q^{u^{(\eta)}}q^{v^{(\eta)}} \Big).\]

It follows that 
\[C_{W_{M}^{'},S_{*}}(q) = \bigg[ \Big(\prod_{u_{1}=1}^\infty \frac{1}{1-q^{u_{1}}} \Big) \Big(\frac{1}{2} \prod_{v_{1}}^\infty \frac{1}{1-q^{v_{1}}} + \frac{1}{2} \prod_{v_{1}}^\infty \frac{1}{1+q^{v_{1}}} \Big) \bigg] \] \\
\[ \times \prod_{\eta \in H_{*} \setminus 1} \bigg[ \Big(\prod_{u_{\eta}=1}^\infty \frac{1}{1-q^{u_{\eta}}} \Big) \Big(\frac{1}{2} \prod_{v_{\eta}}^\infty \frac{1}{1-q^{v_{\eta}}} + \frac{1}{2} \prod_{v_{\eta}}^\infty \frac{1}{1+q^{v_{\eta}}} \Big) \bigg] \] \\
\[= \bigg[ \Big( \frac{1}{2} \prod_{n_{1}=1}^\infty \frac{1}{1-q^{2 n_{1}}} + \frac{1}{2} \prod_{n_{1}=1}^\infty \frac{1}{(1-q^{n_{1}})^{2}} \Big) \bigg] \times  \prod_{\eta \in H_{*} \setminus 1} \bigg[ \Big( \frac{1}{2} \prod_{n_{\eta}=1}^\infty \frac{1}{1-q^{2 n_{\eta}}} + \frac{1}{2} \prod_{n_{\eta}=1}^\infty \frac{1}{(1-q^{n_{\eta}})^{2}} \Big) \bigg] \] \\
\[= \Big( \frac{1}{2} \prod_{k=1}^\infty \frac{1}{1-q^{2 k}} + \frac{1}{2} \prod_{k=1}^\infty \frac{1}{(1-q^{k})^{2}} \Big)^{|H_{*}|}\] \\
\[= \Big(\frac{1}{2} \frac{q^{1/12}}{\eta(z)^{2}} + \frac{1}{2} \frac{q^{1/12}}{\eta(2z)} \Big)^{|H_{*}|}. \]
The equality between the first and seccond line is given in the appendix of \cite{BD}.  The switch of variables to $n_{1}$ and $n_{\eta}$ keeps track that the number of products is indexed by the conjugacy classes.  {\flushright \qed}

\section{Preliminaries}
Here we will recall basic properties of modular forms.  For more information see \cite{ONO1}.

\subsection{Modular forms, operators, and cusps}
A large part of the proof of Theorem \ref{six} involves understanding cusps of congruence subgroups.  A \emph{cusp} of $\Gamma \subset SL_{2}(\Z)$ is an equivalence class of $\Q \cup \{\infty \}$ under the action of $\Gamma$.  We will divide the rest of the section into subsections on integer weight modular forms, half-integral weight modular forms, and then a section on the modularity of eta-quotients.  
\subsubsection{Integer weight modular forms}
In this section let $k$ be an integer which denotes the weight of an integer weight modular form.  For each meromorphic function $f$ on the upper half complex plane $\mathbb{H}$ and each integer $k$, define the slash operator, $\vert_{k}$, by
\[f(z) \vert_{k} \gamma := (\det \gamma)^{k/2} (cz+d)^{-k} f(\gamma z)\]
where $\gamma =  \left( \begin{array}{cc} a & b \\ c & d \end{array} \right) \in GL_{2}^{+}(\R)$.  We say $f$ is a \emph{meromorphic modular form} of weight $k$ on $\Gamma$ if 
\[f \vert_{k} \gamma = f \qquad \emph{for all } \gamma \in \Gamma.\]
We call $f$ a \emph{weakly holomorphic modular form} if its poles are supported at the cusps.  If $f$ is holomorphic at the cusps we say it is a \emph{holomorphic modular form}, and if it vanishes at the cusps we say it is a \emph{cusp form}.  We denote the spaces of weakly holomorphic, holomorphic, and cusp forms with character $\chi$ by $\mathcal{M}_{k}(\Gamma, \chi), M_{k}(\Gamma, \chi),$ and $S_{k}(\Gamma, \chi)$, respectively.  If $\chi$ is a Dirichlet character modulo $N$, then we say that a form $f(z) \in M_{k}(\Gamma_{1}(N))$ (resp. $S_{k}(\Gamma_{1}(N))$ or $\mathcal{M}_{k}(\Gamma_{1}(N))$) has \textit{Nebentypus character} $\chi$ if
\[f \left(\frac{az+b}{cz+d} \right) = \chi(d) (cz+d)^{k} f(z)\]
for all $z \in \mathbb{H}$ and all $\left( \begin{array}{cc} a & b \\ c & d \end{array} \right) \in \Gamma_{0}(N)$, where
\[\Gamma_{0}(N) := \Bigg\{ \left( \begin{array}{cc} a & b \\ c & d \end{array} \right) \in SL_{2}(\Z) : c \equiv 0 \pmod{N} \Bigg\}\]
and \[ \Gamma_{1}(N) := \Bigg\{ \left( \begin{array}{cc} a & b \\ c & d \end{array} \right) \in SL_{2}(\Z) : a \equiv d \equiv 1 \pmod{N}, \rm{and}\  c \equiv 0 \pmod{N} \Bigg\}.\]
The space of such forms is denoted $M_{k}(\Gamma_{0}(N), \chi)$ (resp. $S_{k}(\Gamma_{0}(N), \chi)$ or $\mathcal{M}_{k}(\Gamma_{0}(N), \chi)$).
We will now recall some operators on integer weight modular forms.  Each modular form $f$ has a Fourier expansion $f(z) = \sum a(n) q^n$, where $q := e^{2 \pi i z}$.  If $f(z) = \sum_{n \geq n_{0}} a(n)q^n$, the $U$-operator, $U_{t}$, on $f(z)$ is defined by
\[f(z) \vert U_{t} = \sum_{n \geq n_{0}} a(tn) q^n. \]
Similarly, the $V$-operator, $V_{t}$, is defined by 
\[f(z) \vert V_{t} = \sum_{n \geq n_{0}} a(n) q^{tn}. \]
The following facts can be found in \cite[p. 28]{ONO1}.  Suppose $f(z) \in M_{k}(\Gamma_{0}(N), \chi)$, where $k$ is an integer.
\begin{enumerate}
\item If $t$ is a positive integer, then 
\[f(z) \vert V_{t} \in M_{k}(\Gamma_{0}(Nt), \chi).\]
\item If $t \mid N$, then 
\[f(z) \vert U_{t} \in M_{k}(\Gamma_{0}(N), \chi). \]
\end{enumerate}
Furthermore, if $f(z)$ is a cusp form, then so are $f(z) \vert U_{t}$ and $f(z) \vert V_{t}$.
If $f(z) \in M_{k}(\Gamma_{0}(N), \chi)$, for each prime $p \nmid N$, the integer weight Hecke operator $T_{p,k,\chi}$ preserves the space $ M_{k}(\Gamma_{0}(N), \chi)$ and acts on $f(z) = \sum_{n=0}^\infty a(n) q^{n}$ by 
\[ f(z) \vert T_{p, k, \chi} = \sum_{n=0}^\infty \big( a(pn) + \chi(p) p^{k-1} a(n/p) \big) q^{n}, \]
where $a(n/p) = 0$ if $ p \nmid n$.

\subsubsection{Half-integral weight modular forms}

In this section, let $\lambda$ be an integer, and define $\lambda + \frac{1}{2}$ to be the weight of a half-integral weight modular form.  In order to define the slash operator for half-integral weight modular forms we need to consider the group of pairs $\xi = (\gamma, \phi(z))$ where $\gamma =  \left( \begin{array}{cc} a & b \\ c & d \end{array} \right) \in GL_{2}^{+}(\R)$ and $\phi(z)$ is a complex holomorphic function on $\mathbb{H}$ such that $\phi^{2}(z) = \frac{\pm (cz + d)}{\sqrt{\det(\gamma)}}$.  Define the half-integral weight slash operator by
\[f(z) \vert_{\lambda +\frac{1}{2}} \xi = \phi(z)^{-2\lambda -1} f(\gamma z).\]

For $\gamma =  \left( \begin{array}{cc} a & b \\ c & d \end{array} \right) \in \Gamma_{0}(4)$ and $z \in \mathbb{H}$, define 
\[j(\gamma, z) := \left( \frac{c}{d} \right) \varepsilon_{d}^{-1} \sqrt{cz + d},\]
where $\left( \frac{c}{d} \right)$ is the Jacobi symbol and 
\[\varepsilon_{d} := \begin{cases} 1 & \emph{if } d \equiv 1 \pmod{4} \\
i & \emph{if } d \equiv 3 \pmod{4}.
\end{cases}\]
Set $\tilde{\gamma} := (\gamma, j(\gamma, z))$, and for any congruence subgroup $\Gamma \subset \Gamma_{0}(4)$, let $\tilde{\Gamma} := \{ \tilde{\gamma} : \gamma \in \Gamma \}$.
We say $f$ is a \emph{weakly holomorphic modular form} of weight $ \lambda +\frac{1}{2}$ on $\tilde{\Gamma}$ if it is holomorphic on $\mathbb{H}$, meromorphic at the cusps, and satisfies
\[f \vert_{\lambda +\frac{1}{2}} \tilde{\gamma} = f \]
for all  $\tilde{\gamma} \in \tilde{\Gamma}$.  In the same way as for integer weight, if $f$ is holomorphic at the cusps we say it is a \emph{holomorphic modular form}, and if it vanishes at the cusps we say it is a \emph{cusp form}.  We denote the spaces of weakly holomorphic, holomorphic, and cusp forms with character $\chi$ by $\mathcal{M}_{\lambda +\frac{1}{2}}(\Gamma, \chi), M_{\lambda +\frac{1}{2}}(\Gamma, \chi),$ and $S_{\lambda +\frac{1}{2}}(\Gamma, \chi)$.  If $\chi$ is a Dirichlet character modulo $4N$, then we say $g(z) \in M_{\lambda + \frac{1}{2}}(\Gamma)$ (resp. $S_{\lambda + \frac{1}{2}}(\Gamma)$ or $\mathcal{M}_{\lambda + \frac{1}{2}}(\Gamma)$) has \textit{Nebentypus character} $\chi$ if 
\[g \left(\frac{az+b}{cz+d} \right) = \chi(d) \left(\frac{c}{d} \right)^{2 \lambda + 1} \varepsilon_{d}^{-2 \lambda -1} (cz+d)^{\lambda + \frac{1}{2}} g(z)\]
for all $z \in \mathbb{H}$ and all $\left( \begin{array}{cc} a & b \\ c & d \end{array} \right) \in \Gamma_{0}(4N)$.  The space of such forms is denoted $M_{\lambda + \frac{1}{2}}(\Gamma_{0}(4N), \chi)$ (resp. $S_{\lambda + \frac{1}{2}}(\Gamma_{0}(4N), \chi)$ or $\mathcal{M}_{\lambda + \frac{1}{2}}(\Gamma_{0}(4N), \chi)$).
We will now recall some operators on half-integral weight modular forms.  If $g(z) = \sum_{n \geq n_{0}} b(n)q^n$, the $U$-operator, $U_{t}$, on $g(z)$ is defined by
\[g(z) \vert U_{t} = \sum_{n \geq n_{0}} b(tn) q^n. \]
Similarly, the $V$-operator, $V_{t}$, is defined by 
\[g(z) \vert V_{t} = \sum_{n \geq n_{0}} b(n) q^{tn}. \]
The following can be found in \cite[p. 50]{ONO1}.  Suppose $g(z) \in M_{\lambda +\frac{1}{2}}(\Gamma_{0}(4N), \chi)$.
\begin{enumerate}
\item If $t$ is a positive integer, then 
\[g(z) \vert V_{t} \in M_{\lambda + \frac{1}{2}} \left( \Gamma_{0}(4Nt), \left( \frac{4t}{\bullet} \right) \chi \right).\]
\item If $t \mid N$, then 
\[g(z) \vert U_{t} \in M_{\lambda +\frac{1}{2}} \left( \Gamma_{0}(4N), \left( \frac{4t}{\bullet} \right) \chi \right). \]
\end{enumerate}
Furthermore, if $g(z)$ is a cusp form, then so are $g(z) \vert U_{t}$ and $g(z) \vert V_{t}$.
If $g(z) \in M_{ \lambda + \frac{1}{2}}(\Gamma_{0}(4N), \chi)$, where $\lambda$ is an integer, then for each prime $\ell \nmid 4N$ the half-integral weight Hecke operator preserves the space $ M_{\lambda + \frac{1}{2}}(\Gamma_{0}(4N), \chi)$ and acts on  $g(z) = \sum_{m=0}^\infty b(m) q^m$ by 
\[g(z) \vert T(\ell^{2}, \lambda, \chi) = \sum_{m=0}^\infty \big(b(\ell^{2} m) + \chi^{*}(\ell) \big(\frac{m}{\ell} \big) \ell^{\lambda - 1} b(m) + \chi^{*}(\ell^{2}) \ell^{2 \lambda - 1} b(m/{\ell^{2}}) \big) q^{m}, \]
where $\chi^{*}(m) := (\frac{(-1)^{\lambda}}{m}) \chi(m)$ and $b(m/{\ell^{2}}) = 0$ if $\ell^{2} \nmid m$.

\subsubsection{Modularity of eta-quotients}
\emph{Dedekind's eta-function} is a weight $\frac{1}{2}$ modular form defined as
\[\eta(z) := q^{1/24} \prod_{n=1}^\infty (1 - q^n).\]
An \emph{eta-quotient} is a function $f(z)$ of the form
\[f(z) = \prod_{\delta \mid N} \eta(\delta z)^{r_{\delta}}\]
where $N \geq 1$ and $r_{\delta}$ is an integer.  We have the following useful proposition for eta-quotients.

\begin{Prop}[{\cite[p. 18]{ONO1}}] \label{seven}
If $f(z) = \prod_{\delta \mid N} \eta(\delta z)^{r_{\delta}}$ is an eta-quotient with weight $k = \frac{1}{2} \sum_{\delta \mid N} r_{\delta}$ and the additional properties that 
\[\sum_{\delta \mid N} \delta r_{\delta} \equiv 0 \pmod{24}\]
and 
\[\sum_{\delta \mid N} \frac{N}{\delta} r_{\delta} \equiv 0 \pmod{24},\]
then $f(z)$ satisfies
\[f \left( \frac{az+b}{cz+d} \right) = \chi(d) (cz+d)^{k} f(z)\]
for every $\left( \begin{array}{cc} a & b \\ c & d \end{array} \right) \in \Gamma_{0}(N)$.  Here $\chi(d) := (\frac{(-1)^k s}{d})$ where $s := \prod_{\delta \mid N} \delta^{r_{\delta}}.$

\end{Prop}

\subsection{Theorems of Treneer, Serre, and Ono}

In order to study congruence properties, we turn to a result from Serre on the action of the Hecke operator on integral weight modular forms.

\begin{Prop}[Serre, \cite{Serre2}] \label{eight}

Suppose that $f(z) = \sum_{n=1}^\infty a(n) q^{n} \in S_{k}(\Gamma_{0}(N), \chi)$ has coefficients in $\mathcal{O}_{K}$, the ring of algebraic integers in the number field $K$, and $M$ is a positive integer.  Furthermore, suppose $k >1$.  Then a positive proportion of  the primes $p \equiv -1 \pmod{MN}$ have the property that 
\[f(z) \vert T_{p, k, \chi} \equiv 0 \pmod{M}. \] 

\end{Prop}

There is an analogous proposition for half-integral weight modular forms due to Ono and Ahlgren which is proved using Proposition \ref{eight} and Shimura's correspondence between half-integral weight modular forms and even integer weight modular forms.

\begin{Prop}[{Ahlgren-Ono, \cite[p. 56]{ONO1}}] \label{nine}

Suppose that $g(z) = \sum_{m=1}^\infty b(m) q^{m} \in S_{\lambda + \frac{1}{2}}(\Gamma_{0}(4N), \chi)$ has coefficients in $\mathcal{O}_{K}$, the ring of algebraic integers in the number field $K$, and $M$ is a positive integer.  Furthermore, suppose $\lambda >1$.  Then a positive proportion of  the primes $\ell \equiv -1 \pmod{4MN}$ have the property that 
\[g(z) \vert T(\ell^{2}, \lambda, \chi) \equiv 0 \pmod{M}. \]
\end{Prop}

It is natural to ask for a generalization of Propositions \ref{eight} and \ref{nine} where a Hecke operator for a prime $p$ could simultaneously annihilate a finite set of modular forms.  In order to tackle this problem we will now turn our attention to modular Galois representations.  Let $\overbar{\Q}$ be an algebraic closure of $\Q$, and for each rational prime $\ell$, let $\overbar{\Q}_{\ell}$ be an algebraic closure of $\Q_{\ell}$.  Fix an embedding of $\overbar{\Q}$ into $\overbar{\Q}_{\ell}$.  This fixes  choice of a decomposition group $D_{\ell} = \{ \sigma \in \rm{Gal}(K/\Q) : \sigma(\mathfrak{p}_{\ell, K}) = \mathfrak{p}_{\ell, K} \}$.  Specifically, if $K$ is any finite extension of $\Q$ and $\O_{K}$ is the ring of integers of $K$, then for each $\ell$ this fixes a choice of prime ideal $\mathfrak{p}_{\ell,K}$ of $\O_{K}$ dividing $\ell$.  Let $\mathbb{F}_{\ell,K} = \O_{K}/\mathfrak{p}_{\ell,K}$ be the residue field of $\mathfrak{p}_{\ell, K}$ and let $\vert \cdot \vert_{\ell}$ be an extension to $\overbar{\Q}_{\ell}$ of the usual $\ell$-adic absolute value on $\Q_{\ell}$.

\begin{Thm}[{\cite[p. 42]{ONO1}}] \label{seventeen}
Let $f(z) = \sum_{n=1}^\infty a(n) q^n \in S_{k}(\Gamma_{0}(N), \chi)$ be a newform, and let $K_{f}$ be the field extension obtained by adjoining all of the $a(n)$ and values of $\chi$ to $\Q$.  If $K$ is any finite extension of $\Q$ containing $K_{f}$ and $\ell$ is any prime, then due to work of Eichler, Shimura, Deligne, and Serre there is a continuous semisimple representation
\[ \rho_{f, \ell} : \rm{Gal}(\overbar{\Q}/\Q) \rightarrow GL_{2}(\mathbb{F}_{\ell, K})\]  
for which the following are true:
\begin{enumerate}
\item $\rho_{f, \ell}$ is unramified at all primes $p \nmid N \ell$.
\item $\emph{Tr}(\rho_{f, \ell}(\emph{Frob}_{p})) \equiv a(p) \pmod{\mathfrak{p}_{\ell, K}}$ for all primes $p \nmid N \ell$.
\item $\emph{det}(\rho_{f, \ell}(\emph{Frob}_{p})) \equiv \chi(p) p^{k-1} \pmod{\mathfrak{p}_{\ell, K}}$ for all primes $p \nmid N \ell$.
\item $\emph{det}(\rho_{f, \ell}(c)) = -1$ for any complex conjugation $c$.
\end{enumerate}

\end{Thm}

\begin{rmk}
Let $D_{f, \ell} = G_{f} \cap D_{\ell}$ where $G_{f}$ is the subgroup of $\rm{Gal}(\overbar{\Q}/\Q)$ stabilizing $f$, and let  \[\mathbb{F}_{f, \ell} = \mathbb{F}_{\ell, K}^{D_{f, \ell}} = \{ a \in \mathbb{F}_{\ell, K} : \sigma(a) = a \quad \forall \sigma \in D_{f, \ell} \}.\] 
If $f$ does not have complex multiplication and $\ell$ is sufficiently large, then the image of $\rho_{f, \ell}$ contains a normal subgroup $H_{f}$ conjugate to $SL_{2}(\mathbb{F}_{f, \ell})$.  Essentially, this means that the image of $\rho_{f, \ell}$ is almost always `as large as possible'.

\end{rmk}
Newforms are eigenforms for the Hecke operator $T_{p,k, \chi}$ with eigenvalues given by the $p$th coefficients of the newform.  The fact that the image of $\rho_{f, \ell}$ is large, along with an application of the Chebotarev Density Theorem tells us we can choose the image of $\rm{Frob}_{p}$ to have a trace of zero a positive proportion of the time.  This determines the $p$th coefficient and thus implies Proposition \ref{eight} of Serre.  The following lemma from \cite{OS} extends the idea of these representations having large image and allows us to apply it to sums of modular forms.

\begin{Lemma}[Ono-Skinner, \cite{OS}] \label{eighteen}
Let $f_{1}, f_{2}, ..., f_{v}$ be newforms without complex multiplication, and let $f_{i}(z) = \sum_{n=1}^\infty a_{i}(n)q^n$.  Write $\rho_{f_{i}, \ell} = \rho_{i}$ and $\mathbb{F}_{f_{i}, \ell} = \mathbb{F}_{i}$.  Then
\begin{enumerate}
\item the image of $\rho_{1} \times \cdot \cdot \cdot \times \rho_{v}$ is conjugate to $SL_{2}(\mathbb{F}_{1}) \times \cdot \cdot \cdot \times SL_{2}(\mathbb{F}_{v})$.
\item For each positive integer $d$ and each $w \in \mathbb{F}_{i}$, a positive density of primes $p \equiv 1 \pmod{d}$ satisfies 
\[a_{i}(p) \equiv w \pmod{\mathfrak{p}_{\ell, K}}. \]
\item For each pair of coprime positive integers $r$, $d$, a positive density of primes $p \equiv r \pmod{d}$ satisfies $\vert a_{i}(p) \vert_{\ell} = 1$.
\end{enumerate}

\end{Lemma}

Part $(1)$ of Lemma \ref{eighteen} specifically tells us that, with small adjustments, we can apply Propositions \ref{eight} and \ref{nine} to a finite set of modular forms simultaneously.  This fact is crucial for the proof of Theorem \ref{main}.

A large portion of this paper will apply work of Treneer in \cite{Tre} to $C_{W_{M}^{'},S_{*}^{'}}(q)$.  The main result from \cite{Tre} follows.

\begin{comment} NEED TRENEER????? \end{comment}

\begin{Prop}[Treneer, \cite{Tre}] \label{ten}
Suppose $p$ is an odd prime, and that $k$ and $r$ are integers with $k$ odd.  Let $N$ be a positive integer with $4 \mid N$ and $(N, p) = 1$, and let $\chi$ be a Dirichlet character modulo $N$.  Let $K$ be an algebraic number field with ring of integers $\mathcal{O}_{K}$, and suppose $f(z) = \sum a(n) q^n \in \mathcal{M}_{\frac{k}{2}}(\widetilde{\Gamma_{0}(N)}, \chi) \cap \mathcal{O}_{K}((q))$.  If $r$ is sufficiently large, then for  each positive integer $j$, a positive proportion of primes $Q \equiv -1 \pmod{Np^j}$ have the property that 
\[a(Q^3 p^r n) \equiv 0 \pmod{p^j}\]
for all $n$ coprime to $Qp$.

\end{Prop}

\begin{rmk}
A similar statement is true for integer weight forms which will be apparent in the proof.  This should not be surprising due to the fact that Serre has shown that almost all coefficents of integer weight forms are $0 \pmod{m}$ for any integer $m > 1$.
\end{rmk}

\begin{comment}

\section{Proof of Theorem \ref{five}}

 Given a vector $\textbf{e} := (e_{1}, e_{2}, ..., e_{k}) \in \Z^{k}$ define the \emph{generalized partition function} $p(n)_{\textbf{e}}$ to be the coefficients of the following power series.
\[\sum_{n=0}^\infty p(n)_{\textbf{e}} q^{n} = \prod_{m=1}^{k} \prod_{n=1}^\infty \frac{1}{(1-q^{mn})^{e_{m}}} = \prod_{n=1}^\infty \frac{1}{(1-q^{n})^{e_{1}} \cdot \cdot \cdot (1-q^{kn})^{e_{k}}}.\]
Using a theorem of Tauberian, Cotron, Dicks, and Fleming prove asymptotics for the generalized partition function in \cite{CDF1}.

\begin{Prop}[C-D-F] \label{eleven}
Given a nonzero vector $\textbf{e} := (e_{1}, e_{2}, ..., e_{k}) \in \Z_{\geq 0}^{k}$ let $ d := gcd\{m : e_{m} \neq 0 \}$.  Define $\gamma$ and $\delta$ by 
\[\gamma := \gamma(\textbf{e}) =  \sum_{n=1}^{k/d} e_{dn}, \quad
\delta := \delta(\textbf{e}) =  \sum_{n=1}^{k/d} \frac{e_{dn}}{n}.\]
Then as $n \rightarrow \infty$ we have 
\[p(dn)_{\textbf{e}} \sim \frac{\lambda A^{\frac{1 + \gamma}{4}}}{2 \sqrt{\pi} n^{\frac{3+ \gamma}{4}}} e^{2 \sqrt{An}}, \]
where $\lambda := \prod_{m=1}^{k} \Big( \frac{m}{2 \pi} \Big)^{\frac{e_{dm}}{2}}$ and $A := \frac{\pi^{2} \delta}{6}.$
\end{Prop}

Note that $\sum_{n=0}^\infty c_{M}(n) q^{n} = \frac{1}{2^{M}} \sum_{k = 0}^{M} \binom{M}{k} p(n)_{(2M - 2k, k)} q^{n}$.  Applying the above theorem to this generalized partition function leads to the desired result.  
{\flushright \qed}

\end{comment}

\section{Proof of Theorem \ref{main}}

In order to find congruences for sums of mixed weight modular forms, we must examine where their coefficents overlap.  The following lemma will describe this.

\begin{Lemma} \label{twelve}

Let $f_{i}(z) = \sum_{n=0}^\infty a_{i}(n) q^{n} \in M_{k_{i}}(\Gamma_{0}(N_{i}), \chi_{i})$ and \\ $g_{j}(z) = \sum_{m=0}^\infty b_{j}(m) q^{m} \in M_{\lambda_{j} + \frac{1}{2}}(\Gamma_{0}(4M_{j}), \chi_{j})$. 
 If  \[f_{i}(z) \vert T_{p_{i}} \equiv 0 \pmod{Q} \]
and \[g_{j}(z) \vert T(\ell_{j}^{2}) \equiv 0 \pmod{Q} \] with all $p_{i}$ and $\ell_{j}$ distinct,
then \[a_{1}\left(\prod_{i} p_{i}^{2r_{i} + 1} \prod_{j} \ell_{j}^{4s_{j} + 3} n \right) \equiv \cdot \cdot \cdot \equiv a_{u}\left(\prod_{i} p_{i}^{2r_{i} + 1} \prod_{j} \ell_{j}^{4s_{j} + 3} n \right)\]
\[\equiv b_{1}\left(\prod_{i} p_{i}^{2r_{i} + 1} \prod_{j} \ell_{j}^{4s_{j} + 3} n \right) \equiv \cdot \cdot \cdot \equiv b_{v}\left(\prod_{i} p_{i}^{2r_{i} + 1} \prod_{j} \ell_{j}^{4s_{j} + 3} n \right) \equiv 0 \pmod{Q}, \]
where $gcd(p_{1} \cdot \cdot \cdot p_{u} \ell_{1} \cdot \cdot \cdot \ell_{v}, n) = 1$ and $r_{i}$ and $s_{j}$ are nonnegative integers.  
\end{Lemma}

\begin{rmk}
 if $p_{i} = \ell_{j}$ for some $i$ and $j$, then remove the $p_{i}$ for the congruence to hold.  This is made clear in the following corollary and in the proof of Lemma \ref{twelve}.  

\end{rmk}

\begin{cor}

Let $f_{i}(z)$ and $g_{j}(z)$ be given as above.  If $f_{i}(z) \vert T_{p} \equiv 0 \pmod{Q}$ and $g_{j}(z) \vert T(p^{2}) \equiv 0 \pmod{Q}$ for the same prime $p$, then 
\[a_{1}\left(p^{4t+3} n \right) \equiv \cdot \cdot \cdot \equiv a_{u}\left(p^{4t+3} n \right) \equiv b_{1}\left(p^{4t+3} n \right) \equiv \cdot \cdot \cdot \equiv b_{v}\left(p^{4t+3} n \right) \equiv 0  \pmod{Q}, \]
where $gcd(p, n) = 1$ and $t$ is a nonnegative integer.

\end{cor}

\emph{Proof of Lemma \ref{twelve}:} Assume that $f(z) \vert T_{p} \equiv 0 \pmod{Q}$ and $g(z) \vert T(\ell^{2}) \equiv 0 \pmod{Q}$.  If $p \nmid n$, then $f(z)\vert T_{p, k, \chi} = \sum_{n=0}^\infty a(pn)q^{n} \equiv 0 \pmod{Q}$.  If $p \mid n$ then we can replace $n$ with $p^2 r$ with $r$ and $p$ coprime in $ f(z) \vert T_{p, k, \chi} = \sum_{n=0}^\infty \big( a(pn) + \chi(p) p^{k-1} a(n/p) \big) q^{n}$ to arrive at 
\[\sum_{r=0}^\infty \big( a(p^{3} r) + \chi(p) p^{k-1} a(pr) \big) q^{p^{2}r} \equiv  0 \pmod{Q}. \]  Since we know $a(pr) \equiv  0 \pmod{Q}$ in this case then $a(p^{3}r) \equiv  0 \pmod{Q}$ must be true for $p$ and $r$ coprime.  This process can be repeated to show $a(p^{2t+1}n) \equiv   0 \pmod{Q}$ for $p$ and $n$ coprime. \\
If $\ell \mid m$ and $\ell^{2} \nmid m$, then \[g(z) \vert T(\ell^{2}, \lambda, \chi) = \sum_{m=0}^\infty \big(b(\ell^{2} m) + \chi^{*}(\ell) \big(\frac{m}{\ell} \big) \ell^{\lambda - 1} b(m) + \chi^{*}(\ell^{2}) \ell^{2 \lambda - 1} b(m/{\ell^{2}}) \big) q^{m}\] simplifies to $\sum_{m=0}^\infty b(\ell^{2}m)q^{m}$.  Replacing $m$ with $\ell s$ where $(\ell, s) = 1$ shows that \\
$b(\ell^{3} s) \equiv 0 \pmod{Q}$.  If $\ell^{2} \mid m$, then we can replace $m$ with $\ell^{5} s$ with $\ell$ and $s$ coprime to get
\[\sum_{s=0}^\infty \big(b(\ell^{7} s) + \chi^{*}(\ell^{2}) \ell^{2 \lambda - 1} b(\ell^{3} s) \big) q^{\ell^{5} s} \equiv  0 \pmod{Q}. \]
We know that $b(\ell^{3} s) \equiv  0 \pmod{Q}$ so $b(\ell^{7}s) \equiv  0 \pmod{Q}$ must be true for $\ell$ and $s$ coprime.  This process can be repeated to show $b(\ell^{4t+3}m) \equiv  0 \pmod{Q}$ for $\ell$ and $m$ coprime.
These two observations combined lead to Lemma \ref{twelve}.  {\flushright \qed}

Theorem \ref{main} partly follows from being able to use the above corollary for any set of modular forms.  The congruence in Proposition \ref{ten} comes from Treneer being able to show that any weakly holomorphic modular form has coefficients that are congruent $  \pmod{p^j}$ to the coefficients of a cusp form.  These details are worked out below for the example from section $2$.  Proposition \ref{ten} then follows by applying Proposition \ref{eight} or \ref{nine} to that cusp form.  In order to apply the above corollary, we need to be able to use Propositions \ref{eight} and \ref{nine} on a finite set of modular forms.
Lemma \ref{eighteen} tells us that the images of  the Galois representations of a finite set of modular forms will almost always be simultaneously `as large as possible'.  Due to this fact, we can let $g \in Im( \rho_{1} \times \cdot \cdot \cdot \times \rho_{v})$ be conjugate to 
\[\left( \begin{array}{cc} 0 & 1 \\ -1 & 0 \end{array} \right) \times \cdot \cdot \cdot \times \left( \begin{array}{cc} 0 & 1 \\ -1 & 0 \end{array} \right)\]
and, by the Chebotarev Density Theorem, a positive proportion of primes $Q \equiv -1  \pmod{N p^j}$ satisfy $( \rho_{1} \times \cdot \cdot \cdot \times \rho_{v})(\rm{Frob}_{Q}) = g$.  This allows us to apply Propositions \ref{eight} and \ref{nine} to a finite set of modular forms for the same prime.  Theorem \ref{main} then follows from being able to use techniques from \cite{Tre} for a finite set of weakly holomorphic modular forms to show they are all congruent to cusp forms.  We can then use Propositions \ref{eight} and \ref{nine} on that set of forms to show they are simultaneously annihilated by Hecke operators, so we can apply the above corollary.
{\flushright \qed}

\section{Proof of Theorem \ref{six}}

In this section we will explicitly work out the congruence properties of the conjugacy growth series from section $2$ following the work of Treneer in \cite{Tre}.

Recall that  \[ \sum_{n=0}^\infty \gamma_{W_{M}^{'}, S_{*}^{'}}(n) q^{n} = \left( \frac{1}{2} \frac{q^{1/12}}{\eta(z)^{2}} + \frac{1}{2} \frac{q^{1/12}}{\eta(2z)} \right)^{M}. \]
Define \[F_{M}(z) := \sum_{n=-M}^\infty b_{M}(n)q^{n} = \left( \frac{1}{2} \frac{1}{\eta(12z)^{2}} + \frac{1}{2} \frac{1}{\eta(24z)} \right)^{M} \]
\[= \frac{1}{2^{M}} \sum_{k=0}^{M} \binom{M}{k} \frac{1}{\eta(12z)^{2(M-k)}} \frac{1}{\eta(24z)^{k}} \]
\[= \frac{1}{2^{M}} \sum_{k=0}^{M} \binom{M}{k} F_{M, k}(z), \]
where $F_{M, k}(z) = \sum_{n=-M}^\infty a_{M,k}(n) q^{n} \in M_{\frac{k}{2} - M}(\Gamma_{0}(N_{M, k}))$.
\begin{comment}  We want to use Propositions \ref{eight} and \ref{nine} to show that the $F_{M,k}(z)$ are annihilated by Hecke operators.  However, Propositions \ref{eight} and \ref{nine} require cusp forms.  The following lemma allows us to relate the $F_{M,k}(z)$ to cusp forms in order to use the Propositions.
\end{comment}

\begin{Lemma} \label{thirteen}
Suppose $p$ is an odd prime and $r$ and $N_{M, k}$ are integers with $(N_{M,k}, p) = 1$.  If $r$ is sufficiently large then, for every positive integer $j$ there exists an integer $\beta \geq j -1$ and a cusp form
\[F_{M,k,p,j}(z) \in S_{\frac{k}{2} - M + \frac{p^{\beta}(p^2 -1)}{2}}(\Gamma_{0}(N_{M,k}p^2)) \]
such that 
\[F_{M,k,p,j}(z) \equiv \sum_{p \nmid n} a_{M,k}(p^{r}n) q^{n} (\emph{mod } p^{j}). \]

\end{Lemma}

\emph{Proof of Lemma \ref{thirteen}}: This proof will follow the proof of Proposition $3.1$ in \cite{Tre}.  The plan for this proof is to divide the cusps of $\Gamma_{0}(N_{M,k}p^{2})$ into two groups.  We will pick $r$ large enough so that $F_{M,k}(z) \vert U_{p^{r}} = \sum a_{M,k}(p^{r}n)q^{n}$ is holomorphic at each cusp $\frac{a}{c}$ with $p^{2} \mid c$.  Then we will define \[F_{M,k,r}(z) = F_{M,k}(z) \vert U_{p^{r}} - F_{M,k}(z) \vert U_{p^{r+1}} \vert V_{p} = \sum_{p \nmid n} a_{M,k}(p^{r}n)q^{n}\]
such that it vanishes at these cusps.  
Define the eta-quotients
\[F_{p}(z) = \begin{cases}
\frac{\eta^{p^{2}}(z)}{\eta(p^{2}z)} \in M_{\frac{p^{2} - 1}{2}}(\Gamma_{0}(p)) & p \geq 5 \\
\frac{\eta^{27}(z)}{\eta^{3}(9z)} \in M_{12}(\Gamma_{0}(9)) & p = 3.
\end{cases} \]
$F_{p}(z)$ vanishes at every cusp $\frac{a}{c}$ where $p^{2} \nmid c$ and is $1 \pmod{p}$.  By induction it also clear that $F_{p}(z)^{p^{s-1}} \equiv 1 \pmod{p^{s}}$ for any integer $s$.  Our cusp form will end up being $F_{M,k,r}(z) \cdot F_{p}(z)^{p^{\beta}}$ for some integer $\beta$.
First we must find an explicit description of the Fourier expansion of $F_{M,k}(z) \vert U_{p^{r}}$ at a cusp $\frac{a}{c}$ with $p^2 \nmid c$.  Note that for the remainder of the paper $\gamma$ is used with the slash operator, even when the form being hit with the slash operator is a half-integral weight form.  It is implied that $\gamma$ (and other matrices) should be replaced by $\tilde{\gamma}$ where appropriate.

\begin{Prop} \label{fourteen}

Let $ \gamma = \left( \begin{array}{cc} a & b \\ cp^2 & d \end{array} \right) \in SL_{2}(\Z)$, $\tilde{\gamma} = \left(\left( \begin{array}{cc} a & b \\ cp^2 & d \end{array} \right), \mu \sqrt{cp^2 z + d} \right)$ with $\mu \in \{\pm 1, \pm i \}$ and $ac > 0$.  Then there exists an integer $n_{0} \geq -24M$ and a sequence $\{a_{M,k,0}(n) \}_{n \geq n_{0}}$ such that for each $r \geq 1$,
\[ \big( F_{M,k}(z) \vert U_{p^r} \big) \vert_{\frac{k}{2} - M} \gamma = \sum_{\substack{n=n_{0} \\ n \equiv 0 \pmod{p^r}}}^\infty a_{M,k,0}(n) q_{24 p^r}^{n} \]
where $q_{24 p^r}^{n} = e^{\frac{2 \pi i n z}{24 p^r}}$.

\end{Prop}

\emph{Proof of Proposition \ref{fourteen}}: 
First we note that for any matrix $A \in SL_{2}(\Z)$, $F_{M,k}(z) \vert_{\frac{k}{2} - M} A = \sum_{n=n_{0}}^\infty a_{M,k,0}(n) q_{24}^{n}$ where $n_{0} \geq -24M$.  This can be seen by following Theorem $1$ in \cite{Hon}.  
\[\prod_{i=i}^{j} \left( \eta(t_{i}z) \vert_{\kappa / 2}  \left( \begin{array}{cc} a & b \\ c & d \end{array} \right) \right)^{n_{i}}\]
transforms to 
\[C \prod_{i=i}^{j} \left( \eta(z) \vert_{\kappa / 2}  \left( \begin{array}{cc} \alpha_{i} & \beta_{i} \\ 0 & \delta_{i} \end{array} \right) \right)^{n_{i}}\]
where $\left( \begin{array}{cc} t_{i} & 0 \\ 0 & 1 \end{array} \right) \left( \begin{array}{cc} a & b \\ c & d \end{array} \right) = \left( \begin{array}{cc} a_{i} & b_{i} \\ c_{i} & d_{i} \end{array} \right) \left( \begin{array}{cc} \alpha_{i} & \beta_{i} \\ 0 & \delta_{i} \end{array} \right)$ and $C$ is a constant.  We can see from this that $t_{i}a = \alpha_{i} a_{i}$ and $c= \alpha_{i} c_{i}$, so $\alpha_{i} = (t_{i}, c) \leq t_{i}$.  Taking $t_{1} = 12, n_{1} = -2(M-k), t_{2} = 24,$ and $n_{2} = -k$ we arrive at the conclusion that  $F_{M,k}(z) \vert_{\frac{k}{2} - M} A = \sum_{n=n_{0}}^\infty a_{M,k,0}(n) q_{24}^{n}$ where $n_{0} \geq -24M$.
If we define $\sigma_{v,t} :=   \left( \begin{array}{cc} 1 & v \\ 0 & t \end{array} \right)$, $\widetilde{\sigma_{v,t}} :=  \left( \left( \begin{array}{cc} 1 & v \\ 0 & t \end{array} \right), t^{1/4} \right)$,  then notice that \[F_{M,k}(z) \vert U_{t} = t^{\frac{k-2M}{4} - 1} \sum_{v=0}^{t-1} F_{M,k}(z) \vert_{\frac{k}{2} - M} \sigma_{v,t}. \]  For each $0 \leq v \leq p^{r} - 1$, choose an integer $s_{v} \equiv 0 \pmod{4}$ such that \[s_{v}N_{M,k} \equiv (a + vcp^2)^{-1} (b + vd) \pmod{p^r}\] and set $w_{v} = s_{v} N_{M,k}$.  Note that $w_{v}$ runs through the residue classes of $p^r$ as $v$ does.   Also define 
\begin{align*} \alpha_{v} :&= \left( \begin{array}{cc} a + vcp^2 & \frac{b + vd - aw_{v} - w_{v}vcp^2}{p^r} \\ cp^{r+2} & d -w_{v}cp^2 \end{array} \right), \\
\widetilde{\alpha_{v}} :&= (\alpha_{v}, \mu \sqrt{cp^2 (p^r z - w_{v}) + d})\end{align*}
so that $\sigma_{v,p^r} \gamma = \alpha_{v} \sigma_{w_{v}, p^r}$.  Putting this all together, we have 
\begin{align*} \left( F_{M,k}(z) \vert U_{p^r} \right) \vert_{\frac{k}{2} - M} \gamma &= (p^{r})^{\frac{k-2M}{4} - 1} \sum_{v=0}^{p^{r}-1} F_{M,k}(z) \vert_{\frac{k}{2} - M} \sigma_{v,p^{r}} \gamma \\
&= (p^{r})^{\frac{k-2M}{4} - 1} \sum_{v=0}^{p^{r}-1} F_{M,k}(z) \vert_{\frac{k}{2} - M} \alpha_{v} \sigma_{w_{v},p^{r}}. \end{align*}
In Lemma $3.4$ of \cite{Tre} Treneer shows that $\alpha_{v} \alpha_{0}^{-1} \in \Gamma_{1}(N_{M,k})$.  Beacuse $F_{M,k}(z)$ is invariant under action by $\Gamma_{1}(N_{M,k})$, we now have 
\begin{align*}\left( F_{M,k}(z) \vert U_{p^r} \right) \vert_{\frac{k}{2} - M} \gamma &= (p^{r})^{\frac{k-2M}{4} - 1} \sum_{v=0}^{p^{r}-1} F_{M,k}(z) \vert_{\frac{k}{2} - M} \alpha_{v} \sigma_{w_{v},p^{r}} \\
&= (p^{r})^{\frac{k-2M}{4} - 1} \sum_{v=0}^{p^{r}-1} F_{M,k}(z) \vert_{\frac{k}{2} - M} (\alpha_{v} \alpha_{0}^{-1})^{-1} \alpha_{v} \sigma_{w_{v},p^{r}} \\
&= (p^{r})^{\frac{k-2M}{4} - 1} \sum_{v=0}^{p^{r}-1} F_{M,k}(z) \vert_{\frac{k}{2} - M} \alpha_{0} \sigma_{w_{v},p^{r}}.\end{align*}
Since $\alpha_{0} \in SL_{2}(\Z)$, we have
\begin{align*} \sum_{v=0}^{p^r - 1}  F_{M,k}(z) \vert_{\frac{k}{2} - M} \alpha_{0} \sigma_{w_{v},p^{r}} &= \sum_{v=0}^{p^r - 1} \left( \sum_{n=n_{0}}^\infty a_{M,k,0}(n) q_{24}^{n} \right) \vert_{\frac{k}{2} - M} \sigma_{w_{v}, p^r} \\
&= \sum_{v=0}^{p^r - 1} p^{\frac{-r(k-2M)}{4}} \sum_{n=n_{0}}^\infty a_{M,k,0}(n) e^{\frac{2 \pi i n (z + w_{v})}{24p^r}} \\
&= p^{\frac{-r(k-2M)}{4}} \sum_{n=n_{0}}^\infty a_{M,k,0}(n) q_{24p^r}^{n} \sum_{v=0}^{p^r - 1} e^{\frac{2 \pi i n w_{v}}{24p^r}}.\end{align*}
The numbers $\frac{w_{v}}{24}$ run through the residue classes of $p^r$ as $v$ does, therefore 
\[\sum_{v=0}^{p^r - 1} e^{\frac{2 \pi i n w_{v}}{24p^r}} = \sum_{v=0}^{p^r - 1} e^{\frac{2 \pi i n v}{p^r}}= \begin{cases}
p^r & \emph{if } n \equiv 0 \pmod{p^r} \\
0 & \emph{otherwise},
\end{cases}\]
which gives us 
\[\sum_{v=0}^{p^r - 1}  F_{M,k}(z) \vert_{\frac{k}{2} - M} \alpha_{0} \sigma_{w_{v},p^{r}} = p^{r(1 - \frac{k-2m}{4})} \sum_{\substack{n=n_{0} \\ n \equiv 0 \pmod{p^r}}}^\infty a_{M,k,0}(n) q_{24p^r}^{n}.\]
Therefore,
\[\big( F_{M,k}(z) \vert U_{p^r} \big) \vert_{\frac{k}{2} - M} \gamma = \sum_{\substack{n=n_{0} \\ n \equiv 0 \pmod{p^r}}}^\infty a_{M,k,0}(n) q_{24 p^r}^{n}. \]
{\flushright \qed}

\begin{Prop} \label{fifteen}
Define 
\[F_{M,k,r}(z) := F_{M,k}(z) \vert U_{p^r} - F_{M,k}(z) \vert U_{p^{r+1}} \vert V_{p} \in M_{\frac{k}{2} - M} (\Gamma_{0}(N_{M,k}p^2)). \]
Then for $r$ sufficiently large, $F_{M,k,r}(z)$ vanishes at each cusp $\frac{a}{cp^2}$ of $(\Gamma_{0}(N_{M,k}p^2)$ with $ac > 0$.
\end{Prop}

\emph{Proof of Proposition \ref{fifteen}}:
From Proposition \ref{seven}, we know 
\[\big( F_{M,k}(z) \vert U_{p^r} \big) \vert_{\frac{k}{2} - M} \gamma = \sum_{\substack{n=n_{0} \\ n \equiv 0 \pmod{p^r}}}^\infty a_{M,k,0}(n) q_{24 p^r}^{n} \]
where $n_{0} \geq -24M$.  For $r$ sufficiently large, $-p^r < -24M \leq n_{0}$.  In the Fourier expansion if $a_{M,k,0}(n) \neq 0$, in order for $n \equiv 0 \pmod{p^r}$ to be true, $n \geq 0$ must be true.  Therefore, we have 
\[\big( F_{M,k}(z) \vert U_{p^r} \big) \vert_{\frac{k}{2} - M} \gamma = \sum_{\substack{n \geq 0 \\ n \equiv 0 \pmod{p^r}}} a_{M,k,0}(n) q_{24 p^r}^{n},\]
which shows $F_{M,k}(z) \vert U_{p^r}$ is holomorphic at the cusp $\frac{a}{cp^2}$.  
We will handle the second term in a similar way as the first term in the proof of Proposition \ref{fourteen}.  Define $\tau_{v,t} = \left( \begin{array}{cc} 1 & v/t \\ 0 & 1 \end{array} \right)$, $\widetilde{\tau_{v,t}} = \left( \left( \begin{array}{cc} 1 & v/t \\ 0 & 1 \end{array} \right), 1 \right)$, and note that 
\[F_{M,k}(z) \vert U_{t} \vert V_{t} = t^{-1} \sum_{v=0}^{t-1} F_{M,k}(z) \vert_{\frac{k}{2} - M} \tau_{v,t}. \]
Using this, we see that
\[(F_{M,k}(z) \vert U_{p^r}) \vert U_{p} \vert V_{p} \vert_{\frac{k}{2} - M} \gamma = p^{-1} \sum_{v=0}^{p-1} (F_{M,k}(z) \vert U_{p^r}) \vert_{\frac{k}{2} - M} \tau_{v,p} \gamma.\] 
For each $0 \leq v \leq p-1$, choose $s_{v}^{'} \equiv 0 \pmod{4}$ such that $s_{v}^{'} N_{M,k} \equiv a^{-1}vd \pmod{p}$, and set $w_{v}^{'} = s_{v}^{'} N_{M,k}$.  Define 
\begin{align*} \delta_{v} :&= \left( \begin{array}{cc} 1+aw_{v}^{'}cp + vw_{v}^{'}c^{2}p^{2} & \frac{avd - a^{2}w_{v}^{'}}{p} - acvw_{v}^{'} - bvcp \\ w_{v}^{'}c^{2}p^{3} & 1 - aw_{v}^{'}cp \end{array} \right),\\
\widetilde{\delta_{v}} :&= (\delta_{v}, \sqrt{w_{v}^{'}c^2 p^3 z + 1 - aw_{v}^{'} cp}),\end{align*}
so that $\tau_{v,p} \gamma = \delta_{v} \gamma \tau_{w_{v}^{'},p}$.  In \cite{Tre} Treneer also shows that $\delta_{v} \in \Gamma_{1}(N_{M,k}p)$, so 
\[ (F_{M,k}(z) \vert U_{p^r}) \vert_{\frac{k}{2} - M} \tau_{v,p} \gamma =  (F_{M,k}(z) \vert U_{p^r}) \vert_{\frac{k}{2} - M} \delta_{v} \gamma \tau_{w_{v}^{'}, p} =  (F_{M,k}(z) \vert U_{p^r}) \vert_{\frac{k}{2} - M} \gamma \tau_{w_{v}^{'}, p}. \]
Following the same method as in the proof of Proposition \ref{fourteen}, we can write 
\begin{align*} (F_{M,k}(z) \vert U_{p^r}) \vert U_{p} \vert V_{p} \vert_{\frac{k}{2} - M} \gamma &= p^{-1} \sum_{v=0}^{p-1} (F_{M,k}(z) \vert U_{p^r}) \vert_{\frac{k}{2} - M} \gamma \tau_{w_{v}^{'}, p} \\
&= p^{-1} \sum_{v=0}^{p-1} \left( \sum_{\substack{n \geq 0 \\ n \equiv 0 \pmod{p^r}}} a_{M,k,0}(n) q_{24p^r}^{n} \right) \vert_{\frac{k}{2} -M} \tau_{w_{v}^{'},p} \\
&= p^{-1} \sum_{v=0}^{p-1}  \sum_{\substack{n \geq 0 \\ n \equiv 0 \pmod{p^r}}} a_{M,k,0}(n) \rm{exp}\left({\frac{2 \pi i n (z+ \frac{w_{v}^{'}}{p})}{24p^r}} \right) \\
&= \sum_{\substack{n \geq 0 \\ n \equiv 0 \pmod{p^r}}} a_{M,k,0}(n) q_{24p^r}^{n} \sum_{v=0}^{p-1} \rm{exp} \left({\frac{2 \pi i n w_{v}^{'}}{24p^{r+1}}} \right) \\
&= \sum_{\substack{n \geq 0 \\ n \equiv 0 \pmod{p^r}}} a_{M,k,0}(n) q_{24p^r}^{n} \sum_{v=0}^{p-1} \rm{exp} \left({\frac{2 \pi i w_{v}^{'}}{24p}\left(\frac{n}{p^r}\right)} \right).\end{align*}
The numbers $\frac{w_{v}^{'}}{24}$ run through the residue classes modulo $p$ as $v$ does, so
\[\sum_{v=0}^{p-1} \rm{exp} \left({\frac{2 \pi i w_{v}^{'}}{24p}\left(\frac{n}{p^r}\right)} \right) = \sum_{v=0}^{p-1} \rm{exp} \left( \frac{2 \pi i v}{p} \left(\frac{n}{p^r} \right) \right) = \begin{cases}
p & \emph{if } p \mid \frac{n}{p^r} \\
0 & otherwise. 
\end{cases}\]
Putting everything together gives us 
\[(F_{M,k}(z) \vert U_{p^r}) \vert U_{p} \vert V_{p} \vert_{\frac{k}{2} - M} \gamma = \sum_{\substack{n \geq 0 \\ n \equiv 0 \pmod{p^{r+1}}}} a_{M,k,0}(n) q_{24p^r}^{n}.\]
To finish the proof of Proposition \ref{fifteen}, we have 
\[F_{M,k,r}(z) \vert_{\frac{k}{2}-M} \gamma =  \sum_{\substack{n \geq 0 \\ n \equiv 0 \pmod{p^r}}} a_{M,k,0}(n) q_{24p^r}^{n} -  \sum_{\substack{n \geq 0 \\ n \equiv 0 \pmod{p^{r+1}}}} a_{M,k,0}(n) q_{24p^r}^{n}.\]
The constant terms (which may be $0$) of each expansion cancel, so $F_{M,k,r}(z)$ vanishes at the cusp $\frac{a}{cp^2}$.
{\flushright \qed}

Before discussing the cusp $\frac{a}{c}$ where $p^2 \nmid c$, notice that \[F_{M,k,r}(z) = \sum_{n=1}^\infty a_{M,k}(p^r n) q^{n} - \sum_{n=1}^\infty a_{M,k}(p^{r+1} n) q^{pn} = \sum_{\substack{n=1 \\ p \nmid n}}^\infty a_{M,k}(p^r n) q^{n}.\] 

Recall  the eta-quotients
\[F_{p}(z) = \begin{cases}
\frac{\eta^{p^{2}}(z)}{\eta(p^{2}z)} \in M_{\frac{p^{2} - 1}{2}}(\Gamma_{0}(p)) & p \geq 5 \\
\frac{\eta^{27}(z)}{\eta^{3}(9z)} \in M_{12}(\Gamma_{0}(9)) & p = 3,
\end{cases} \]
and recall that $F_{p}(z)$ vanishes at every cusp $\frac{a}{c}$ where $p^{2} \nmid c$.  The forms $F_{p}(z)$ are $1 \pmod{p}$, and by induction it is easy to show $F_{p}(z)^{p^{s-1}} \equiv 1 \pmod{p^{s}}$ for any integer $s$.  Let $r$ be sufficiently large, and fix $j$.  If $\beta \geq j - 1$ is sufficiently large, then 
\[F_{M,k,p,j}(z) := F_{M,k,r}(z) \cdot F_{p}(z)^{p^{\beta}} \equiv F_{M,k,r}(z) \pmod{p^j}\]
vanishes at all cusps $\frac{a}{c}$ of $\Gamma_{0}(N_{M,k} p^2)$ where $p^2 \nmid c$.  By Proposition \ref{fifteen}, $F_{M,k,r,p,j}(z)$ also vanishes at the cusps $\frac{a}{c}$ where $p^2 \mid c$, so 
\[F_{M,k,p,j}(z) \in S_{\frac{k}{2} - M + \frac{p^{\beta}(p^2 - 1)}{2}}(\Gamma_{0}(N_{M,k} p^2)).\] 
As seen above,
\[F_{M,k,p,j}(z) \equiv F_{M,k,r}(z) \equiv \sum_{\substack{n=1 \\ p \nmid n}}^\infty a_{M,k}(p^r n) q^{n} \pmod{p^j}.\]
This completes the proof of Lemma \ref{thirteen}. {\flushright \qed}

\begin{Lemma} \label{sixteen}

\begin{enumerate}
\item If $k$ is even, then $F_{M,k,p,j}(z)$ is an integer weight cusp form, so for a positive proportion of primes $Q \equiv -1 \pmod{N_{M,k} p^j}$, we have 
\[a_{M,k}(Q^{2t+1} p^{r} n) \equiv 0 \pmod{p^j}\]
for all nonegative integers $t$, and $n$ coprime to $Qp$.
\item If $k$ is odd, then $F_{M,k,p,j}(z)$ is a half-integral weight cusp form, so for a positive proportion of primes $Q \equiv -1 \pmod{N_{M,k} p^j}$, we have
\[a_{M,k}(Q^{4t+3} p^r n)  \equiv 0 \pmod{p^j}\]
for all nonegative integers $t$, and $n$ coprime to $Qp$.
\end{enumerate}

\end{Lemma}

\emph{Proof of Lemma \ref{sixteen}}:
If $k$ is even (resp. odd), then $F_{M,k,p,j}(z)$ is an integral (resp. half-integral) weight cusp form.  Thus, by Proposition \ref{eight} (resp. Proposition \ref{nine}), for a positive proportion of primes $Q \equiv -1 \pmod{N_{M,k} p^j}$, we have $F_{M,k,p,j}(z) \vert T_{Q} \equiv 0 \pmod{p^j}$ (resp. $F_{M,k,p,j}(z) \vert T(Q^2) \equiv 0 \pmod{p^j}$).  If we let $F_{M,k,p,j}(z) = \sum_{n=1}^\infty c_{M,k}(n)q^{n}$, then by Lemma \ref{twelve} we have $c_{M,k}(Q^{2t+1}n) \equiv  0 \pmod{p^j}$  (resp. $c_{M,k}(Q^{4t+3}n) \equiv  0 \pmod{p^j}$) for any nonnegative integer $t$ and $Q$ and $n$ coprime. The rest of the proof follows from the fact that $c_{M,k}(n) \equiv a_{M,k}(p^r n) \pmod{p^j}$. {\flushright \qed}

We will now refer back to part ($1$)  of Lemma \ref{eighteen}.  Using this and Chebotarev's Density Theorem, we are able to apply Proposition \ref{eight} or \ref{nine} simultaneously to each term in a sum of modular forms, which in turn allows us to apply Lemma \ref{twelve} to our entire sum at the same time instead of piece by piece as in Lemma \ref{sixteen}.  As in Theorem \ref{main}, if we have a sum of modular forms $f_{i}$ of mixed weights and level $N_{i}$, we can replace the level $N$ in Proposition \ref{eight} or \ref{nine} with the smallest $N'$ such that each $N_{i}$ divides $N'$.

\begin{comment}
\begin{Thm} \label{nineteen}
Suppose $p \geq 5$ be prime and let $j \geq 1$.  Then if $r$ is sufficiently large for a positive proportion of primes $Q \equiv -1 (\emph{mod } 576p^{j})$, we have that 
\[c_{M} \left(\frac{Q^{4t+3}p^{r}n + M}{12} \right) \equiv 0 (\emph{mod } p^{j})\]
for all $n$ coprime to $Qp$, and for all nonnegative integers $t$.

\end{Thm}
\end{comment}
We will now complete the proof of Theorem \ref{six}.  Recall that 
\[ \sum_{n=0}^\infty \gamma_{W_{M}^{'}, S_{*}^{'}}(n) q^{n} = \left( \frac{1}{2} \frac{q^{1/12}}{\eta(z)^{2}} + \frac{1}{2} \frac{q^{1/12}}{\eta(2z)} \right)^{M} \]
and 
\begin{align*} F_{M}(z) := \sum_{n=-M}^\infty b_{M}(n)q^{n} &= \left( \frac{1}{2} \frac{1}{\eta(12z)^{2}} + \frac{1}{2} \frac{1}{\eta(24z)} \right)^{M} \\
&= \frac{1}{2^{M}} \sum_{k=0}^{M} \binom{M}{k} \frac{1}{\eta(12z)^{2(M-k)}} \frac{1}{\eta(24z)^{k}} \\
&= \frac{1}{2^{M}} \sum_{k=0}^{M} \binom{M}{k} F_{M, k}(z), \end{align*}
so $F_{M, k}(z) = \sum_{n=-M}^\infty a_{M,k}(n) q^{n} \in M_{\frac{k}{2} - M}(\Gamma_{0}(N_{M, k}))$.  Note also that 
\[\sum_{n=-M}^\infty b_{M}(n)q^{n} = \sum_{k=0}^{M} \sum_{n=-M}^\infty \binom{M}{k} a_{M, k}(n)q^{n}.\]
In this sum there will be a form of level $576$ and all of the other forms will have level dividing $576$.  Clearly, the sum will be a mix of integer weight and half-integral weight modular forms.
From Lemma \ref{sixteen} we know that for a positive proportion of primes $Q \equiv -1 \pmod{N_{M,k} p^j}$, we have 
\[a_{M,k}(Q^{2t+1} p^{r} n) \equiv 0 \pmod{p^j}\] for $k$ even and 
\[a_{M,k}(Q^{4t+3} p^r n)  \equiv 0 \pmod{p^j}\] for $k$ odd,
for all nonegative integers $t$, and $n$ coprime to $Qp$.  Theorem \ref{seventeen} and Lemma \ref{eighteen} together imply that Lemma \ref{sixteen} can be applied to each $F_{M,k}(z)$ simultaneously for a positive proportion of primes $Q \equiv -1 \pmod{576p^j}$, so $a_{M,k}(Q^{4t+3}p^r n) \equiv 0 \pmod{p^j}$ for each $a_{M,k}(n)$.  Since the congruence holds for each part of the sum, we also have 
\[b_{M}(Q^{4t+3}p^r n) \equiv 0 \pmod{p^j}\] for a positive proportion of primes $Q  \equiv -1 \pmod{576p^j}$.  Because $b_{M}(n) = \gamma_{W_{M}^{'}, S_{*}^{'}} \left(\frac{n+M}{12} \right)$, we have 
\[\gamma_{W_{M}^{'}, S_{*}^{'}} \left(\frac{Q^{4t+3}p^{r}n + M}{12} \right) \equiv 0 \pmod{p^{j}}.\]
{\flushright \qed}

\section{Details of the example}

At the end of Section $1$ the following example was given:
\begin{ex}
We find that
\[\gamma_{W_{2}^{'}, S_{*}^{'}}\left(\frac{7 n + 2}{12}\right) \equiv 0 \pmod{7}\]
whenever $n \not\equiv 10 \pmod{24}$.  Moreover, the above congruence is true when $ n = 24t + 10$ and $t \equiv 2, 4, 5,$ or $6 \pmod{7}$.  
\end{ex}
We will now give the details of this example.  Define
\[\sum_{n=0}^\infty a_{1}(n) q^{n} = \prod_{n=1}^\infty \frac{1}{(1-q^{n})^4},\]
\[\sum_{n=0}^\infty a_{2}(n) q^{n} = \prod_{n=1}^\infty \frac{1}{(1-q^{n})^2} \prod_{n=1}^\infty \frac{1}{(1-q^{2n})^{2}},\] and 
\[\sum_{n=0}^\infty a_{3}(n) q^{n} = \prod_{n=1}^\infty \frac{1}{(1-q^{2n})^2}.\]
It is clear that 
\[\sum_{n=0}^\infty \gamma_{W_{2}^{'}, S_{*}^{'}} (n) q^n = \frac{1}{4} \sum_{n=0}^\infty \left( a_{1}(n) + a_{2}(n) + a_{3}(n) \right) q^n.\]
By adapting Theorem $6$ from \cite{Ono2}, we have
\[\sum a_{1} \left( \frac{pn +2}{12} \right) q^n \equiv \frac{\Delta^{\frac{p^2 - 1}{6}}(z) \vert U_{p} \vert V_{12}}{\eta^{4p}(12z)} \pmod{p},\]
\[\sum a_{2} \left( \frac{pn +2}{12} \right) q^n \equiv \frac{\Delta^{\frac{p^2 - 1}{12}}(z) \Delta^{\frac{p^2 - 1}{24}}(2z) \vert U_{p} \vert V_{12}}{\eta^{2p}(12z) \eta^{p}(24z)} \pmod{p},\] and 
\[\sum a_{3} \left( \frac{pn +2}{12} \right) q^n \equiv \frac{\Delta^{\frac{p^2 - 1}{12}}(2z) \vert U_{p} \vert V_{12}}{\eta^{2p}(24z)} \pmod{p},\]
where $\Delta(z) := \eta^{24}(z)$ is \textit{Ramanujan's tau function}.  Using a theorem of Sturm in \cite{St}, one can verify with a finite computation that 
\[\Delta^{8}(z) \vert U_{7} \equiv 0 \pmod{7},\]
\[\Delta^{4}(z) \Delta^{2}(2z) \vert U_{7} \equiv 0 \pmod{7},\] and
\[\Delta^{4}(2z) \vert U_{7} \equiv 3 \Delta(2z) \pmod{7}.\]
From this it is clear that 
\[\sum \gamma_{W_{2}^{'}, S_{*}^{'}}\left(\frac{7 n + 2}{12}\right) q^n \equiv \frac{3 \Delta(24z)}{\eta^{14}(24z)} \equiv 3 \eta^{10}(24z) \pmod{7}.\]
In \cite{Serre}, Serre showed that $\eta^{10}$ is lacunary, so one should expect a lot of congruences.  In \cite{WL}, Locus and the author use this fact to prove the congruences necessary to complete the example.

\end{document}